\documentclass[a4paper,11pt,twoside]{article}
\usepackage[hscale=0.7,centering]{geometry}
\usepackage{times}
\usepackage{amssymb,amsmath}
\usepackage{epsf}
\usepackage{times,bm,mathrsfs}  
\usepackage{amsmath}
\usepackage{amssymb}
\usepackage{amsfonts}
\usepackage{mathrsfs}
\usepackage{bbm}
\usepackage{color}
\usepackage{graphicx}
\usepackage{amscd}
\usepackage{epsfig}

\allowdisplaybreaks

\definecolor{Red}{rgb}{1,0,0}
\definecolor{Blue}{rgb}{0,0,1}
\definecolor{Olive}{rgb}{0.41,0.55,0.13}
\definecolor{Green}{rgb}{0,1,0}
\definecolor{MGreen}{rgb}{0,0.8,0}
\definecolor{DGreen}{rgb}{0,0.55,0}
\definecolor{Yellow}{rgb}{1,1,0}
\definecolor{Cyan}{rgb}{0,1,1}
\definecolor{Magenta}{rgb}{1,0,1}
\definecolor{Orange}{rgb}{1,.5,0}
\definecolor{Violet}{rgb}{.5,0,.5}
\definecolor{Purple}{rgb}{.75,0,.25}
\definecolor{Brown}{rgb}{.75,.5,.25}
\definecolor{Grey}{rgb}{.5,.5,.5}
\definecolor{Black}{rgb}{0,0,0}

\def\path{{\tt path}}

\newcommand{\fcal}{\mathcal{F}}

\newcommand{\tccal}{$\mathcal{C}$}

\newcommand{\ga}{\alpha}

\newcommand{\bM}{\mathbf M}
\newcommand{\bT}{\mathbf T}
\newcommand{\gb}{\beta}

\newcommand{\gd}{\delta}
\newcommand{\gD}{\Delta}

\newcommand{\ind}{\mathbbm{1}}

\newcommand{\Wn}{W^{(t)} }
\newcommand{\Xn}{X^{(t)} }
\newcommand{\Yn}{Y^{(t)} }
\newcommand{\Zn}{Z^{(t)} }

\newcommand{\bdm}{\begin{displaymath}}
\newcommand{\edm}{\end{displaymath}}
\newcommand{\bea}{\begin{eqnarray*}}
\newcommand{\eea}{\end{eqnarray*}}
\newcommand{\bean}{\begin{eqnarray}}
\newcommand{\eean}{\end{eqnarray}}

\newcommand{\prob}{\mathbb{P}}

\newcommand{\E}{\mathbb{E}}
\newcommand{\I}{\mathbb{I}}

\newcommand{\mom}[1]{{\bar #1}}

\newcommand{\bfxx}{\mathbf{X}}

\newtheorem{theorem}{Theorem}

\newtheorem{corollary}{Corollary}
\newtheorem{definition}{Definition}

\newtheorem{lemma}{Lemma}

\newtheorem{remark}{Remark}

\title{Branching Process approach for 2-SAT thresholds
}

\author{Elchanan Mossel\thanks{Dept. of Statistics,
    U.C. Berkeley. Supported by an Alfred Sloan fellowship in
    Mathemamatics, NSF grants DMS 0528488 and DMS 0548249 (Career),
and ONR grant N0014-07-1-05-06.}
\and
Arnab Sen \thanks{Dept. of Statistics, U.C. Berkeley. Supported by NSF grants DMS 0528488 and DMS 0548249 and ONR grant N0014-07-1-05-06.}
}
\begin{document}

\maketitle

\thispagestyle{empty}

\begin{abstract}
 It is well known that,  as $n$ tends to infinity, the probability of satisfiability for a random 2-SAT formula on $n$ variables, where each clause occurs independently with probability $\ga/2n$, exhibits a sharp threshold at $\ga=1$.
 We  study a more general 2-SAT model in which each clause
occurs independently but with probability $\ga_i/2n$ where
 $i \in \{0,1,2 \}$ is the number of positive literals in that clause.
 We  generalize  branching process arguments by Verhoeven(99) to determine the satisfiability  threshold
 for this model in terms of the maximum eigenvalue of the branching
 matrix.
\end{abstract}

\bigskip


\clearpage

\section{Introduction}

\subsection{Background}
 The $k$-satisfiability (in short k-SAT)
problem is a canonical constraint satisfaction problem in
theoretical computer science. A $k$-SAT formula is a conjunction
of $m$ clauses, each of which is a disjunction of length $k$
chosen from $n$ boolean variables and their negations. Given a
$k$-SAT formula a natural question is to find an assignment of $n$
variables which satisfies the formula. The decision version of the
problem is to determine whether there exists an assignment
satisfying the formula.

From the computational complexity perspective, the problem is well
understood. The problem is NP-hard for $k \geq 3$ \cite{cook71}
and linear time solvable for $k=2$ \cite{apt79}. Much recent
interest was devoted to the understanding of random $k$-SAT
formulas where each clause is chosen independently with the same
probability and the expected number of clauses in the formula is
$\alpha n$. This problem lies in the intersection of three
different subjects --- statistical physics, discrete mathematics
and complexity theory.

In statistical physics, the notion of `phase transition' refers to
a situation where systems undergo some abrupt behavioral change
depending on some external control parameter such as temperature.
In the context of random $k$-SAT formulas the natural parameter is
the density of the formula $\alpha$, i.e., the ratio between the
number of clauses to the number of variables. Much recent research
is devoted to understanding the critical densities for random
$k$-SAT problems. The most important one is the critical density
for satisfiability, i.e., the threshold at which a formula becomes
from satisfiable with high probability to unsatisfiable with high
probability \cite{Mezard06, yuval03}. Other thresholds involve the
geometry of the solution space and the performance of various
algorithms (e.g.\ see \cite{ach06, mossel07} and references
therein).

The problem of $2$-SAT is more amenable to analysis than $k$-SAT
for $k \geq 3$. This is closely related to the fact that $2$-SAT
can be solved in linear time and to a clear graph theoretic
criteria which is equivalent to satisfiability. The threshold for
$2$-SAT is well known to be $\alpha=1$ ( see \cite{Goe92, Goe96,
CR92, vega92} ) and detailed information on the scaling window is
given in \cite{wilson}. For $k$-SAT for $k \geq 3$ there are
various bounds and conjectures on the critical threshold for
satisfiability but the thresholds are not known rigorously
\cite{friedgut, yuval03, ach00, mezard02}.

In the current paper we establish the threshold of a more general
$2$-SAT model where the probability of having a clause in the
formula depends on the number of positive and negative variables
in the clause. Our proof is based on branching process arguments.
Branching process techniques have been used before to study the
standard $2$-SAT formulas in the unsatisfiability regime ($\alpha
>1$), for example see \cite{verh99}. We generalize the arguments
given in \cite{verh99} in the $2$-type branching process set-up to
analyze the general $2$-SAT model. Our main contribution is in
demonstrating that branching process arguments extend to a
multi-type setup. A well accepted idea in studying random graphs
and constraint satisfaction problems is that since ``local''
structure of the problems is tree-like, processes defined on trees
play a key role in analyzing the problems. The classical example
is the threshold for the existence of a ``giant'' component in
random graphs where branching processes play a key role in the
proof (see e.g \cite{jansonbook, bollobasbook}). Some more recent
examples include \cite{Aldous:01, Weitz:06, MoWeWo:08}.

A seemingly closely related work is by Cooper et. al.
\cite{cooper02} where the threshold for random $2$-SAT with given
literal degree distribution is derived. We note that the two
papers are incomparable, since ours is stated in terms of the
distribution of clauses of different types while in
\cite{cooper02} the model is stated in terms of the degrees of the
literals. While the distribution of clauses determines the
distribution of the degrees of the literals, the converse does not
hold. For example, a random $2$-SAT formula with $2n$
positive-negative clauses has the same literal degree distribution
as a uniform random formula with $2n$ clauses. It is obvious that
while the former is always satisfiable, the latter is not
satisfiable with high probability.

 \subsection{Definitions and statements of main results}

Let  $x_1, x_2, \ldots, x_n$ be $n$ boolean \emph{variables} .
Calling  the negation of $x_i$ as $\bar x_i$, these $n$ boolean
variables give us $2n$ \emph{literals} $\{ x_1, x_2, \ldots , x_n,
\mom x_1, \mom x_2, \ldots, \mom x_n \}$. The two literals $x_i$
and $\bar x_i$ are called complementary to each other ($x_i =1$
iff $ \bar x_i = 0$) with the convention that ${ \bar{\bar x}}_i=
x_i$. We will call the $n$ literals $x_i, i=1, 2, \ldots, n$
\emph{positive literals} and their complementary literals $\mom
x_i, i=1, 2, \ldots, n$ \emph{negative literals}.

Given a literal $u$,  $var(u)$ denotes the corresponding variable, the notation naturally extending to a set of literals $S$ by
$var(S) = \{ var(u): u \in S\} $.  Two literals $u$ and $v$ are said to \emph{ strongly distinct} if $u \ne v$ and $ u \ne \bar v$,
 equivalently, if $var(u) \ne var(v)$. A \emph{2-clause} (which we will call simply a ``clause'' later) is a disjunction $C = x \vee y $ of  two strongly distinct literals. In this paper we will not allow $ x \vee x$ or $ x \vee \bar x$ as valid clauses.  A  2-SAT formula is a conjunction $F = C_1 \wedge C_2 \wedge \ldots \wedge C_m$ of 2-clauses $C_1,C_2, \ldots, C_m$. Let $\mathcal C = \{ C_1,C_2, \ldots, C_m \} $ be  the collection of clauses corresponding to $F$.

As usual in the boolean algebra, $0$ stands for the logical value FALSE, and $1$ stands for the logical value TRUE. A 2-SAT formula $ F = F(x_1, x_2, \ldots, x_n)$ is said to be \emph{satisfiable} if there exists a  truth assignment  $ \eta = (\eta_1, \eta_2, \ldots, \eta_n) \in \{0,1 \}^n$ such that $F(\eta_1, \eta_2, \ldots, \eta_n)  =1 $. The formula $F$ is called \emph{SAT} if the formula $F$ is satisfiable and $F$ is called \emph{UNSAT} otherwise.

In the standard model for a random 2-SAT formula we choose each of
the  possible $ 4 { n \choose 2}$ clauses independently with
probability $ \ga/2n$. In this paper we will study a more general
model in which a 2-SAT formula $F$ consists of  a random subset
$\mathcal C$ of  clauses such that each clause appears in
$\mathcal C$ independently and a clause having $i$ positive
literals is present in the formula with probability $\ga_i/2n, i
=0,1,2$ for some constants $\ga_i \ge 0$.  Of course, taking
$\ga_0 = \ga_1 = \ga_2 = \ga$ we retrieve the standard model.

Let $\bM$ be the branching matrix given by
\begin{equation}\label{eq:branchmatrix}
        \bM =
            \frac{1}{2}\left[ \begin{array} {cc}
             \ga_1 & \ga_0\\
              \ga_2 & \ga_1
            \end{array} \right].
            \end{equation}
Note that though $\bM$ is not symmetric in general, its eigenvalues are all real and given by
 $\frac{1}{2}(\ga_1 \pm \sqrt{\ga_0 \ga_2})$. Let  $\rho = \frac{1}{2}(\ga_1 + \sqrt{\ga_0 \ga_2})$ denote the largest eigenvalue of $M$. We show that $\rho$  is
 a crucial parameter for satisfiability. In particular, our main result is the following theorem which establishes that the generalized 2-SAT to model undergoes a phase transition from satisfiability to unsatisfiability at $\rho =1$.

 \begin{theorem}\label {main:thm}
Let $F$ be random 2-SAT formula under generalized model with parameter $\rho$. \\
(a) If $\rho < 1$ or $\ga_0 \ga_2  =0 $ then $F$ is satisfiable with probability tending to one as $n \to \infty$.\\
(b) If $\rho > 1$ and $\ga_0 \ga_2  >0$ then $F$ is unsatisfiable with probability tending to one as $n \to \infty$.
\end{theorem}

\begin{remark}
It is easy to see and well known that the satisfiability
threshold for 2-SAT remains the same for variants of the model where the set of
2-clauses contains also clauses of the form $x \vee y$ where $x$ and $y$ may not be
strongly distinct. Similarly, the threshold remains the same
if $n \alpha $ clauses are chosen uniformly at random instead of
choosing each clause independently with probability $\ga/2n$
(See the appendix A of \cite{wilson}). The same reasoning applies to the more general 2-SAT model considered here.
\end{remark}

\section{2-SAT and the implication digraph}
We will exploit the standard representation of a 2-SAT formula as a directed graph
(see \cite{wilson} for example), called the {\em implication digraph} associated
 with the 2-SAT formula. This graph has $2n$ vertices, labelled by the $2n$ literals.
If the clause $ (u \vee v) $ is present in the 2-SAT formula then
we draw  the two directed edges  $\mom u \to  v$ and $\mom v \to
u$. The directed edges can be thought of as logical implications
since if there is a directed edge from $u \to v$ and $u =1$, then
for the formula to be satisfiable, it is necessary to have $v=1$.

By a \emph{directed path} (from $u$ to $v$), we mean a sequence of vertices $ u_0 = u, u_1, u_2, \ldots, u_k= v$ such that
there is a directed edge from $ u_i \to u_{i+1}$ for all $i =0,1,\ldots, k-1$. The length of this directed path is $k$.
 A \emph{contradictory cycle} is a union of
 two (not necessarily disjoint) directed paths - one starts from a literal $u$ and ends at its compliment $\bar u$ , the other starts from $\bar u$ and ends at $u$.

The following lemma connects the concept of
 satisfiability of the 2-SAT problem
 to the existence of contradictory cycle in the implication digraph. For a proof, see \cite{wilson}.

\begin{lemma} \label{lem:contradic_cyc}
 A 2-SAT formula is satisfiable iff its implication diagraph contains no contradictory cycle.
\end{lemma}

\section{ Proof of Theorem \ref{main:thm} part (a)}
The proof has some resemblance with the first moment arguments given in Chv\'atal and Reed \cite{CR92}. In fact our `hooked chain' is same as what they called `bicycle'.  The extension to the more general case considered here uses a recursive argument which allows us to deal with the multi-parameter general model.

\begin{definition} Suppose there exists strongly distinct literals
$y_1, y_2, \ldots, y_s$ and $u, v \in$ $ \{y_1, y_2, \ldots, y_s,
 \mom y_1$,  $\mom y_2, \ldots, \mom y_s \}$ such that $(\mom u \vee y_1), (\mom y_1 \vee y_2), \ldots, (\mom y_{s-1} \vee y_s) , (\mom y_s \vee v) \in$ \tccal \  or equivalently, there exists $u \to y_1\to y_2 \to \cdots \to y_s \to v$
 in the implication digraph corresponding to the 2-SAT formula.
We call this sequence of literals a \emph{hooked chain} (of length $s+1$).
\end{definition}

\begin{lemma}\label{lem:hook}
If a 2-SAT formula is unsatisfiable then its implication digraph contains a
 hooked chain of length $\ge 3$.
 \end{lemma}

\noindent \textbf{Proof.} Suppose a 2-SAT formula is unsatisfiable. By lemma \ref{lem:contradic_cyc}, we have  a  contradictory cycle in the implication digraph, say $ u_0  \to u_1 \to u_2 \to \cdots \to u_l = \mom u_0 \to u_{l+1} \to u_{l+2} \to \cdots \to u_{k} = u_0$.  The cycle has at least one directed path from a literal to its complement. Choosing one that minimizes the length we get  an implication chain formed by a sequence of the literals $u_h  \to u_{(h+1) \text{mod }  k} \to u_{(h+2)  \text{mod }  k} \to \cdots \to u_{ (h+t) \text{mod } k}   = \mom u_h $  so that $ u_{(h+1) \text{mod }  k}, u_{(h+2)  \text{mod }  k}, \cdots, u_{ (h +t) \text{mod } k} $ are strongly distinct. Find the largest $s \ge t$ such that  $ u_{(h+1) \text{mod }  k}, u_{(h+2)  \text{mod }  k}, \cdots, u_{ (h +s) \text{mod } k} $ are strongly distinct. Let $v$ be the element pointed to by $u_{ (h +s) \text{mod } k}$ in the cycle. Then  clearly $u_h \to  u_{(h+1) \text{mod }  k} \to u_{(h+2)  \text{mod }  k} \to \cdots \to u_{ (h +s) \text{mod } k}  \to v$ is a hooked chain of length $s+1$. Since there can be no edge between a literal $w$ and its complement $\bar w$, we must have $t \ge 2$ and therefore, $s \ge 2$.  \hfill $\square$

\begin{lemma}
\begin{equation}\label{eq:unsat}
\prob(F \mbox{ is unsatisfiable}) \le  C \sum_{s=2}^{n} \frac{
(2s)^2}{n} [T_{s-1}^+ + T_{s-1}^-]
\end{equation}
where $T_{s-1}^+ (\mbox{or } T_{s-1}^- )$  is  the expected number
of directed paths  of length $(s-1)$ started from $x_1$ (or $\mom
x_1$) consisting of strongly distinct literals  with $T_{0}^+  =
T_{0}^- = 1$ and $C= [\max (\ga_0, \ga_1, \ga_2)]^2$.
\end{lemma}
\noindent \textbf{Proof.}  Let $H_{s} $ be the number of hooked chain of length $(s+1)$ in the implication digraph for 2-SAT  and $\Gamma_{s} $  be the  number of directed paths of $s$ strongly distinct literals in the same digraph. From Lemma \ref{lem:hook},
\begin{align}
\prob(F \mbox{ is unsatisfiable}) &\le  \prob(\exists \mbox{ a hooked chain of length } (s+1)  \mbox{ for some } s \ge 2)  \label{in:1} \\
& \le \sum_{s=2}^{n} \E(H_{s} ) \label{in:2}  \\
& \le  C\sum_{s=2}^{n} \frac{ (2s)^2}{(2n)^2}\E(\Gamma_{s} ) \label{in:3} \\
& = C \sum_{s=2}^{n} \frac{ s^2}{n} [T_{s-1}^+ + T_{s-1}^-]
\label{step4}
\end{align}
 Step \eqref{in:2} follows from simple union bound and Markov inequality. Inequality \eqref{in:3} can be explained as following:

$H_{s}  =  \sum^*_{ y_1, \cdots, y_s }\sum_{  u, v \in  \{y_1, \cdots, y_s,
 \mom y_1, \cdots, \mom y_s \} } \text{I} ( (\bar u \vee y_1), ( \bar y_1\vee y_2), \cdots,  ( \bar y_s \vee v) \in \mathcal C  ) .$

Here $\sum^* $ means the sum is taken over all possible set of strongly distinct literals of size $s$. Observe that for a fixed choice of strongly distinct literals there are $ 2s$ choices for each $u$ and $v$ and each clause
 occurs with probability at most $\max(\ga_0, \ga_1, \ga_2) /2n$. Now  taking expectation and using independence between clauses, we have
 \begin{align*}
 \E( H_s) &=  \sum^*_{ y_1, \cdots, y_s }\sum_{  u, v \in  \{y_1, \cdots, y_s,
 \mom y_1, \cdots, \mom y_s \} } \prob ( ( \bar y_1\vee y_2), \cdots,  ( \bar y_{s-1} \vee y_s) \in \mathcal C  ) \prob ( (\bar u \vee y_1) \in \mathcal C  ) \prob (  ( \bar y_s \vee v) \in \mathcal C  ) \\
 &\le \frac{C (2s)^2}{(2n)^2}  \sum^*_{ y_1, \cdots, y_s } \prob ( ( \bar y_1\vee y_2), ( \bar y_2 \vee y_3), \cdots,  ( \bar y_{s-1} \vee y_s) \in \mathcal C  )\\
 &=  \frac{C (2s)^2}{(2n)^2}  \E (\Gamma_s).
 \end{align*}

Noting that  the quantities $T_{s-1}^+ $ or $ T_{s-1}^- $ do not depend on $x_1$, step \eqref{step4} follows.

\hfill $\square$

\begin{lemma}\label{lem:iter}
Write $\bT_{k} = (T_{k}^+ , T_{k}^-)^T $. Then
\begin{equation}\label{eq:main}
\bT_{s-1} \le  \bM^{s-1} \ind
\end{equation}
where $\bM$ is defined in $\eqref{eq:branchmatrix}$ and $\ind = (1,1)^T$ .
\end{lemma}

\noindent \textbf{Proof.} For a  literal $u$ strongly distinct from $x_1$, let  $J_u$ denote the number of directed paths  of length $(s-2)$ starting from $u$ consisting of strongly distinct literals and not involving the variable $ x_1 $. Then

\begin{align*}
T_{s-1}^+ &= \sum_{\substack{u: u \text{ literals }\\ var(u) \ne var(x_1)}}
\E [ J_u  \times \I ((\mom x_1 \vee u) \in \text{\tccal} )] \\
&= \sum_{\substack{u: u \text{ literals }\\ var(u) \ne var(x_1)}}
\E [ J_u] \times
  \prob ((\mom x_1 \vee u) \in \text{\tccal} ).
\end{align*}

The last step follows from the independence of clauses.

 Simple coupling argument yields  $\E [J_u]   \le T_{s-2}^{+} \text{ or } T_{s-2}^{-}$ depending whether
  the literal  $u$ is positive or negative.
Combining the above facts, we get the following recursive inequality,
\begin{align}
T_{s-1}^+ &\le n T_{s-2}^+  \prob(\mom x_1 \vee x_2)+ n T_{s-2}^-  \prob(\mom x_1 \vee \mom x_2) \notag \\
&\le  \frac{\alpha_1}{2} T_{s-2}^+ + \frac{\alpha_0}{2} T_{s-2}^-
\end{align}
and similarly,
\begin{equation}
T_{s-1}^- \le   \frac{\alpha_2}{2} T_{s-2}^+ + \frac{\alpha_1}{2} T_{s-2}^-
\end{equation}

 Now the above two equations can be written in a more compact way as follows
\begin{equation}\label{eq:iterate}
\bT_{s-1} \le  \bM \bT_{s-2}
\end{equation}

Iterating \eqref{eq:iterate},  we get $\bT_{s-1} \le \bM^{s-1} \bT_{0}= \bM^{s-1}\ind  $.  \hfill $\square$

\vspace{ 3mm}

\noindent {\bf Proof of theorem \ref{main:thm} part (a).} We are now ready to finish the proof of  part (a) of Theorem \ref{main:thm}. If  $\ga_0 \ga_2 =0$ then either all zero or all one assignment always satisfies it. So, take $\ga_0 \ga_2 >0$. Then $\bM$ is semisimple (i.e.\ similar to a diagonal matrix).

By lemma \ref{lem:iter}, $T_{s-1,n}^+  + T_{s-1,n}^- \le \ind^T \bM^{s-1}\ind \le B \rho^{s-1} $ for some constant $B$. The last inequality holds since we assume $\bM$ is semisimple.
Plugging it in \eqref{eq:unsat}, we finally have
\begin{align}
\prob(F \mbox{ is unsatisfiable}) &\le \frac{K}{n} \sum_{s=2}^n s^2\rho^{s-1} \text { for some constant }K>0 \notag \\
 &\le O(n^{-1}) \quad \text{ since } \rho < 1.
\end{align}
\hfill $\square$

\section{The Exploration Process }

Observe that when $\rho > 1$, we need to find a contradictory
cycle in the implication digraph of the random 2-SAT formula with
high probability. In order to prove this, we will show that
starting from any fixed vertex there is a constant probability
that  it  \emph{implies} a large number of literals in the
digraph, meaning that there are directed paths to a large number
of vertices from the fixed vertex. To achieve this, we explore the
digraph dynamically starting from a fixed literal $x$ under
certain rules and keep track of variables that are implied by $x$
at each step. We call this the {\em exploration process} which is
defined next.

\noindent \textbf{Definition and Notations.}
Given a realization of the 2-SAT formula and an arbitrarily fixed literal $x$,  we will consider an exploration process in its implication digraph starting from $x$.
\begin{enumerate}

\item [$\bullet$] This process describes the evolution of two sets of literals which will be called  the \emph{ exposed set} and  the \emph{ active set}.

\item [$\bullet$] A literal is said to be \emph{alive} in a particular step of the process if it is strongly distinct from those in the exposed set and from those in the active set at that step.

\item [$\bullet$]  We maintain two stacks for the literals in active set, one for positive literals and another for negative ones.

\item [$\bullet$]  At each step  we pop-up a literal (call it the \emph{current} literal)  from one of the two stacks of active set, depending on some event to be described later and  \emph{ expose} it. It means that we look for all the  literals that are alive at that time and  to which there is a directed edge from the current literal.
\item [$\bullet$] We then put those new literals in the stacks of the active set ( positive or negative) in some predetermined order and throw the current literal in the exposed set.

\item [$\bullet$]  We go on repeating this procedure until the stack of the active literals becomes empty  and the process stops.

\item [$\bullet$]  Mathematically, let $E_t$ and $A_t$ denote  respectively the set of exposed and  active set of literals at time $t$. Also let $U_t = \{ var(u) : var(u) \not \in var(E_t) \cup var(A_t) \}$ be the set of alive variables at time $t$. Set $ E_0 = \emptyset , A_0 =\{ x \} $.
If $A_t$ is
 non-empty and the literal $ l \in A_t$ is exposed from the stack, then we have the following updates at time $t+1$,
 $$A_{t+1}= (A_t \setminus \{l\}) \cup \{ u : u  \mbox{ literal s.t. }
    var(u) \in U_t  \ \mbox{and the clause} \ (\mom l \vee u) \ \mbox{is present} \},\quad E_{t+1}  = E_t  \cup \{ l\} .$$
     If $A_t$ is empty, then so is $A_{t+1}$ and $E_{t+1}$ will be same as $E_t$.
\end{enumerate}

 Note that during the evolution of the process, each clause is examined only
 once. Also every  literal in $\bigcup_t (A_t \cup E_t)$  can be reached from $x$ via a directed path (consisting of  strongly distinct literals).

For a subset $S$ of literals, we can partition it as $S= S^+ \cup
S^-$ where $S^+$ (resp. $S^-$) is the set of all positive ( resp.
negative) literals of $S$. Let $u_t, a_t^+, a_t^-$ be the
shorthand for $|U_t|$, $|A_t^+|$ and $|A_t^-|$ where $|\cdot|$
means the size of a set. Set $a_t := |A_t| =  a_t^+ + a_t^-$.

\vspace{5mm} \noindent \textbf{Distribution of the process.}  The
stochastic description of the  evolution of the process $(u_t,
a_t^+, a_t^-)$, $0 \le t \le n$ for a random 2-SAT formula on $n$
variables can be summarized in the next lemma whose proof is
immediate.

\begin{lemma} \label{lem:dist}
Define a triangular array of independent Bernoulli random variables.
$$ \Wn_i \sim \mbox{Ber}( \ga_0/2n ); \  \  \Xn_i, \Yn_i \sim \mbox{Ber}( \ga_1/2n  ); \ \  \Zn_i \sim \mbox{Ber}( \ga_2/2n ), \quad 1 \le i \le n, 0 \le t \le n.$$
Let $A_t \ne \emptyset$. Given $ H(t)$, the history up to $t$ and that the current literal  at time $t$ is  positive, we have
\[ u_t - u_{t+1} \stackrel{d}{=}  \sum_{i=1}^{u_t} \{(\Wn_i +\Xn_i) \wedge 1 \}=   \mbox{Bin} (u_t, \ga_0/2n+ \ga_1/2n -\ga_0 \ga_1/4n^2 )\]
\[a_{t+1}^+ - a_t^+   \stackrel{d}{=} -1 + \sum_{i=1}^{u_t} \Xn_i = -1 + \mbox{Bin} ( u_t,  \ga_1/2n), \quad
a_{t+1}^- -   a_t^-   \stackrel{d}{=}   \sum_{i=1}^{u_t} \Wn_i =  -1 +  \mbox{Bin} ( u_t, \ga_0/2n ).
 \]
Similarly, given $H(t)$ and conditional on the event that the current literal  at time $t$  is  negative, we have
\[  u_{t} - u_{t +1} \stackrel{d}{=}  \sum_{i=1}^{u_t} \{ (\Yn_i +\Zn_i) \wedge 1\}=  \mbox{Bin} ( u_t, \ga_1/2n+\ga_2/2n -\ga_1\ga_2/4n^2) \]
\[a_{t+1}^+ - a_t^+   \stackrel{d}{=} -1 +  \sum_{i=1}^{u_t} \Zn_i  =  -1 +  \mbox{Bin}( u_t, \ga_2/2n), \quad
a_{t+1}^- - a_t^-   \stackrel{d}{=}   \sum_{i=1}^{u_t} \Yn_i  =  -1 + \mbox{Bin}( u_t,  \ga_1/2n).
 \]
 \end{lemma}

 \begin{definition}\label{def:T}
 For the rest of the paper, we fix $T = [\sqrt n]$. Let $\tau := \sup \{t  \le T: u_t  \ge u_0 - 2 \ga T \}$ where $\ga = \max( \ga_0, \ga_1, \ga_2)$. In words, $\tau$ is the
last time before $T$  such that the decrease in the number of
unexposed variables is at most $2\ga T$.

   We define a \emph{round} of the exploration process as follows. We fix a subset $S$ of variables of size $N  \ge (1-\gd/2)n$ for some small $\gd > 0$ such that $(1-\delta) \rho >1$ and   a starting literal, say $x \in S$ .   We first run the exploration process from $x$ on the implication digraph restricted to  literals from $S$  up to time $\tau$. If $\tau < T $, we stop. Otherwise, we delete all the variables in $ var( (E_T \cup A_T) \setminus var(x)) $ from $S$ to get a new set of variables $S' \subseteq S$. By the definition of $\tau$,  $| S'| \ge N-2\ga T $. Now we again run another independent exploration process starting from $\mom x$ up to time $\tau$ but on the digraph restricted to literals from $S'$. We again have the stopping rule of $\tau < T$.
 \end{definition}

 \begin{lemma}\label{lem:ind}
The (random) set of clauses examined during the evolution of an exploration process up to time $T$  is  disjoint
 from the set $\{ (\bar u \vee \bar v): u, v \in A_T \} $. Further, the clauses examined during the evolution of
 the second exploration process of a round  are distinct from the clauses in the set  $\{ (\bar u \vee \bar v): u, v \in A_T \} \cup \mathcal D$ where  $A_T$ (resp. $\mathcal D$)  is  the  active set at time $T$ (resp. the set of clauses examined during the evolution ) of the first  exploration process.
 \end{lemma}

\noindent \textbf{Proof:} The first statement of the lemma follows from the easy observation that if $ u \in A_T$, then, from the very definition of the exploration process,  the clause $ (\bar u \vee w)$ is not examined up to time $T$ for all literals $ w$ such that $\bar w \not \in  E_T$.

For the second statement, note that any problematic clause should
include literal $x$ or $\bar x$ ($x= $ starting vertex) . Now all
the clauses involving $var(x)$ which are examined during the first
exploration process of a round  must have the form $ (\bar x \vee
y)$ where literal $y$ is such that $ var(y) \ne var(x)$.  But the
clauses involving $var(x)$ scanned by second exploration process
of the round  are all of the form $( x \vee y)$ where literal $y$
is such that $ var(y) \ne var(x)$ and there can be no clause from
the set  $\{ (\bar u \vee \bar v): u, v \in A_T \}$ which contains
either $x$ or $\bar x$. \hfill $\square$

\begin{corollary}\label{cor:ind}
 Given whether each of the clauses in $\{ (\bar u \vee \bar v): u, v \in A_T \} \cup \mathcal D$ is present in $\mathcal C$ or not , the  distribution of the evolution the second exploration process only depends on number of  the variables  with which the second process starts.
\end{corollary}

\subsection{ Proof of the Theorem \label{thm:main}  part (b) when $ \ga_i$' s are all equal}
Before tackling the general situation we pause for a moment to
give a quick sketch, after \cite{verh99},  of unsatisfiability
part of the phase transition
 for the standard 2-SAT model.  This will serve as a prelude to the proof for the general
 case.

 Let $\alpha = \ga_0 = \ga_1 = \ga_2 >1$. Then $\rho = \ga$.
  In this special case, we slightly modify our exploration process by demanding that we will always
  choose the current literal from the set of active literals uniformly. Thus at each time $t \ge 1$,
  given its size $a_t$, $A_t$ is uniformly random over all the literals except $x$ and $\mom x$, the starting vertex and its complement.

Since the probabilities for the clauses to be present are all equal, the distribution of the exploration process $( u_t, a_t)$ simplifies. Given  $H(t)$, the history up to time $t$ and $a_t > 0 $,
\begin{equation}
 u_t - u_{t+1} \stackrel{d}{=}    \mbox{Bin} (u_t, 2p_{n} -p_{n}^2 ), \quad
a_{t+1} - a_t  \stackrel{d}{=}    -1 + \mbox{Bin}( 2u_t,  p_{n}), \text{ where } p_n = \ga/2n.
\end{equation}
Note that each of the random variables $(u_t - u_{t+1})$ is
stochastically dominated by $\mbox{Bin}(n, 2p_n)$ which has mean $
\alpha$. Thus using concentration of the Binomial distribution  it
is easy to see  that  for time $T = [\sqrt n ]$, the event $\{
\tau < T \} = \{\sum_{t=1}^{T-1}(u_{t-1} - u_{t})
> 2 \alpha T\} $ occurs with probability at most $ A \exp( -cT) $
where $c > 0$. See Lemma \ref{lem:tau} in Section
\ref{subs:twotype} for a proof of a more general fact.

Let $\delta >0$ be as given in Definition \ref{def:T}. If $ u_0
\ge (1-\gd/2)n$, then $ \{ \tau =T \} \subseteq \{ u_T \ge (1-\gd
)n \}$. When $\tau = T $, the process $ \{a_t, t\ge 1\}$ behaves
like a random walk with positive drift on nonnegative integers
with $0$ as the absorbing state and hence
\begin{equation}\label{eq:lb}
\exists C>0 \text{ such that }  \prob( a_T \ge C T, \tau =T ) \ge
\zeta  \text{ for some constant } \zeta > 0, \text{ independent of
$n$}.
\end{equation}

If  both $u$ and $\mom u$ are in $A_T$ for some literal $u$, we have a directed path from the starting vertex to its complement
using the literals in $E_T \cup A_T$. Else, each pair of literals $u, v \in A_T$ are strongly distinct.  There are edges
  $u \to \mom v$ and  $v \to \mom u$ in the digraph if the clause $(\mom u \vee \mom v)$ is present in the formula.  If at least one of these $C T \choose 2$ many clauses is present, then  we again  have a directed path from the starting vertex to its complement using the literals in $E_T \cup A_T$. Let $D$ be the event that there exists  a directed path from  the starting vertex of the exploration process to its complement  in $E_T \cup A_T$. Therefore, by Lemma \ref{lem:ind},
\[ \prob(D \big | \  a_T \ge CT,
 \tau= T ) \ge 1- (1 -\ga/2n)^{CT \choose 2} \ge p, \ \text{where } p > 0 \text{ is a constant, independent of $n$}.\]
which implies that
\[ \prob(D , \tau = T ) \ge  \prob( D \big | \  a_T \ge CT,
 \tau= T )  \prob( a_T \ge C T, \tau =T )  \ge  p \xi > 0. \]

Now, by Corollary \ref{cor:ind}, we can say that after a round, the probability that there is no termination and
there exists a  contradictory cycle in the variables visited during the round  is at least $ p^2 \zeta^2 $.

We continue with another round of exploration process in the deleted graph containing only unvisited variables.
We repeat this process until a round stops due to the stopping rule. If  each of the successive rounds does not
terminate, we have $\Theta(\sqrt n)$ rounds of the exploration processes before the event $\{ u_t < (1- \gd/2)n \} $
occurs. It is easy to see the clauses examined in different rounds of exploration process are all distinct and hence the rounds are independent.

Thus, the probability that we get no contradictory cycle in all the rounds is at most
\begin{align*}
&\quad \  \ \prob( \text{no contradictory cycle and no round stops}) + \prob( \text{ one of the rounds stops})\\
&\le (1-p^2 \zeta^2)^{\Theta(\sqrt n) } +  (1- (1 - A \exp( -cT))^{2\Theta(\sqrt n)} )\\
&\le \text{ Const} \times  \exp {(-B \sqrt n)}  \quad \text{for
some  } B > 0.
\end{align*}
\hfill $\square$

\begin{remark}\label{re:simpleproof} Instead of taking $T =[\sqrt n] $ as in the proof  if we choose
$T  = \Theta(n)$ suitably, then \eqref{eq:lb} yields that $a_T \ge
\Omega(n)$ with probability at least $r$ with some $r>0$. Thus for
any literal $y \ne x, \mom x$,  we get $\prob( y \in A_T ) \ge p$
for some $p>0$ and for all $n$ large enough. So, the probability
that there is a directed path from $x_1$ to $x_2$ is at least $p$.
The same holds true for directed path from $x_2$ to $x_1$. These
are monotonic events. So, by the FKG inequality \cite{fkg}, they
occur simultaneously with probability  greater than or equal to
$p^2$. Thus the chance that there exists a directed path from
$x_1$ to $\mom x_1$ is at least $p^2$. Again applying FKG, we have
a contradictory cycle with probability at least $p^4$. Now
appealing to Friedgut's theorem on sharp threshold
\cite{friedgut}, we can conclude that the formula is UNSAT with
probability tending to one as $n \to \infty$.

\end{remark}

\subsection{Associated 2-type Branching
Process}\label{subs:twotype}

Now we are back to the general case. Given an exploration process on a subgraph of the implication digraph consisting of  $N= \Theta(n)$ many variables starting from any fixed literal, our goal is to couple it with a suitable 2-type supercritical branching process up to time $T = [\sqrt n]$ on a set of high probability. Assume $N \ge (1- \gd/2)n $ where $\delta >0$ is such that $(1-\delta) \rho >1$.

Thus on the set where $\{ u_t  \ge N - 2 \ga T, \ \ \forall t \le T\}  $, for large enough $n$, Bin$(u_t, \ga_i/2n)$  stochastically dominates \ Bin$((1-\delta)n , \ga_i/2n)$ for all time $t \le T$.   Next we are going to prove that this event happens with high probability.

\begin{lemma}\label{lem:tau}
Let $T, \delta$ be as above and  $\ga = \max( \ga_0, \ga_1,
\ga_2)$. Then $\prob( \tau <T ) \le 2\exp(-\ga T/2)$. Therefore,
$\prob(  \ u_t \ge (1 - \gd)n  \ \ \forall t \le T ) \ge 1 -
2\exp(-\ga T/2) $ \ \ for sufficiently large  $n$.
\end{lemma}

\noindent \textbf{Proof.} Since $u_t $ is decreasing in $t$ and $ N - 2 \ga T  \ge (1 - \gd)n$ for sufficiently large  $n$, it is enough to prove that
\[\prob( \tau <T ) \le  2\exp(-\ga T/2). \]
 Note that $u_0 = N$. Clearly, $u_{t-1} - u_{t}$  is conditionally independent of   $u_0 -u_1, u_1- u_2, \ldots, u_{t-2} - u_{t-1}$  given $u_{t-1}$ and the type the current literal at time $t$ and is stochastically dominated by Bin$(2N, \ga/2n)$ irrespective of the conditioning event . Therefore, the distribution of  $u_0 - u_T $ is stochastically dominated by Bin$( 2NT, \ga/2n)$.
By Bernstein's inequality,
\[ \prob(\text{Bin}( 2NT, \ga/2n) \ge 2 \ga T) \le 2 \exp(-\ga T/2) .\]
Therefore,
\[ \prob( \tau <T )  =  \prob ( \ u_T  < N - 2 \ga T )   =  \prob( u_0 - u_T > 2\ga T) \le 2 \exp(-\ga T/2).\]
\hfill $\square$

\begin{lemma}\label{lem:brp}
There exist bounded distributions $F_i$ with mean $m_i,  0 \le i \le 2$ such that

\begin{enumerate}
\item For all sufficiently large $n$
$$\prob(\text{\emph{Bin}}(n(1-\delta) , \ga_i/2n) =k) \ge \prob( X =k | X \sim F_i) \ \forall k \ge 1.$$
\item If $\bM_0$ is the branching matrix given by
\begin{equation}\label{eq:branchmatrix0}
        \bM_0 =
            \left[ \begin{array} {cc}
             m_1 & m_0\\
              m_2 & m_1
            \end{array} \right],
            \end{equation}
 then $\rho_0 > 1$ where $\rho_0$ is the maximum eigenvalue of $\bM_0$.
\end{enumerate}

\end{lemma}

\noindent \textbf{Proof.} Fix some $\gb \in (0,1)$ so that $(1-\delta) (1-\gb)\rho >1$.  Let $ \gamma_i  =(1-\delta)  \ga_i/2$ for $i=0,1,2$. Find $c$ large enough so that
\[ \sum_{k=1}^{c} k (1-\gb/2) \exp(- \gamma_i) \gamma_i^k / k!  >  (1-\gb) \gamma_i \quad i=0,1,2 \]

For each $0\le i \le 2$, let us now define a truncated (and reweighed) Poisson distribution which takes the value $k$ with probability $(1-\gb/2) \exp(- \gamma_i) \gamma_i^k / k!$ for $1 \le k \le c$ and $0$ otherwise. Call this distribution $F_i$. By the choice of $c$, its mean $m_i$ is greater than $(1-\gb) \gamma_i$. Poissonian convergence says Bin$(n(1-\delta) , \ga_i/2n) \stackrel{L^1}{\to}  \mbox{Poisson}(\gamma_i) $ and conclusion $1$ of the lemma follows.

Also, $\rho_0  = m_1 + \sqrt{ m_0 m_2} \ge (1-\delta) (1-\gb)\rho > 1$.
\hfill $\square$

\begin{definition} \label{process} Let us have a supercritical 2-type branching process, which we will call an F-branching process, with offspring distributions as follows
\begin{equation}\label{eq:offdist}
1\mbox{st type} \to (1\mbox{st type}, 2\mbox{nd type}) :   \underbrace{(F_1, F_0)}_{indep}  \quad 2\mbox{nd type}\to (1\mbox{st type},  2\mbox{nd type}) :    \underbrace{(F_2, F_1)}_{indep}
\end{equation}

Next we define  a new process $\bfxx(t) = ( X_1(t), X_2(t) )$ by sequentially traversing the Galton-Watson tree of the F-Branching process. We fix a suitable order among the types of nodes of the tree and moreover we always prefer to visit a node of type I to a node of type II.  Then we traverse the tree sequentially and at each step we expand the tree by including all the children of the node we visit. Let us denote number of unvisited or unexplored children of type $i$ in the tree traversed up to time $t$ by $X_i(t)$.
\end{definition}

\begin{lemma}\label{dominate}
There exists a  coupling such that $(a_t^+, a_t^-) \ge \bfxx(t)$ for all $t \le \tau$ and $n$ large enough.
\end{lemma}

\noindent \textbf{Proof.} Fix $n$ sufficiently large. If the starting vertex of exploration process  is of positive type, we initiate the branching process with one individual of type I. Similarly for the other case. We run in parallel the exploration process where the choice of the type of the current literal at time $t$ depends on the  type of the visited node at time $t$. It can be done because, if $ t \le  \tau$,  we can always simultaneously choose our random variables in such a way (by Lemmas \ref {lem:dist} and \ref{lem:brp})  that for every step, the number  of active literals generated of  each type  is no less than the number of unvisited nodes of the corresponding type in the tree grown up to that step. If $\tau < t \le T$ or if we have no unvisited child left in the tree then we choose the current literal from the active set in some fixed predetermined procedure.
\hfill $\square$

Next we are going to find a lower bound on the total number of unvisited children after $T$ steps of the above process.

\begin{lemma}\label{lem:sc}
Suppose $\bfxx(t)$ be as in Definition \ref{process}  with $\bfxx(0) = (1,0)$  or $ (0,1)$. Then $ \exists C> 0, \eta >0$ such that \ $\prob( X_1(T) + X_2(T) \ge C T) \ge \eta$.
\end{lemma}

\noindent \textbf{Proof.} Though a proof of the above lemma  can be found implicitly  in \cite{ks67}, we here present it for sake of completeness. Recall that  the F-branching process is supercritical as $\rho_0$, the maximum eigenvalue of $\bM_0$, is strictly greater than $1$.

If we assume $\ga_1 >0$, trivially, this process is positive
regular and non-singular. Thus, by a well-known result (see
\cite{harris}) on supercritical multitype branching process, its
extinction probability is given by $\mathbf 0 \le \mathbf q = (
q_1, q_2) < \mathbf 1$ where $q_i$ is the probability that the
process becomes extinct starting with one object of type $i$.

Note if  $\ga_1 =0$, we no longer have the positive regularity.  In that case, though the above theorem can not be directly applied, we can argue as follows to get the same conclusion. If the process starts with only one individual of type $i$,  the corresponding branching process can be viewed as a single type supercritical branching process ( made of the individuals of type $i$ only) if we observe the process only at the even number of steps. So, the probability that it eventually dies out, which is nothing but $q_i$, is strictly less than one.

Let us denote $e_1 =(1, 0) , e_2 =(0,1)$. Instead of looking at $\bfxx(t)$ which has $(0,0)$ as an absorbing state, we will consider a new chain
$\hat \bfxx(t)$ starting from $ \hat \bfxx(0) = \bfxx(0)$ which is supported on entire $\mathbb Z^2$. Given $\hat \bfxx(0), \hat \bfxx(1), \cdots, \hat \bfxx(t)$, define

\begin{equation}
 \hat \bfxx(t+1) \stackrel{d}{=} \left \{ \begin{matrix} \hat \bfxx(t) -e_1 + (F_1, F_0) & \text{if } \hat X_1(t) >0 \\   \hat \bfxx(t) - e_2 + (F_2, F_1) & \text{o.w.}  \end{matrix} \right.
\end{equation}
We can couple $\bfxx(t)$ and $ \hat \bfxx(t)$  together so that  $ \hat \bfxx(t) = \bfxx(t)$ until $ \bfxx(t)$ reaches $(0,0)$.

Let $(a,b) $ be a normalized eigenvector of $\bM_0$ corresponding to eigenvalue $\rho_0$ so that $ a^2 +b^2 = 1$.  Since  $\ga_0, \ga_2 >0$ we have both $a> 0$ and $b>0$. Let $ Z(t) := a X_1(t) + b X_2(t)$ and $\hat Z(t) := a\hat X_1(t) + b\hat X_2(t)$. Let $ \fcal_t := \sigma(\hat \bfxx(0), \hat \bfxx(1), \cdots, \hat \bfxx(t))$ and $\gD \hat Z(t) :=\hat Z(t+1)  - \hat Z(t)$. Since $F_i$'s are bounded, so are $\gD \hat Z(t)$'s. Then

\begin{equation}\label{eq:submart1}
\E(\gD \hat Z(t) | \fcal_t) = \left\{  \begin{array}{lc} (\rho_0-1)a  & \mbox{ if }  \hat X_1(t) >0 \\  (\rho_0-1)b & \mbox{ o.w.} \end{array}\right. \ge \mu := (\rho_0-1) \min(a,b) > 0
\end{equation}

Now we have,
\begin{align*}
\prob ( \hat Z(T) \le \mu T/2) & \le \prob \left ( \hat Z(T) - \sum_{i=0}^{T-1}\E ( \gD \hat Z(i)|\fcal_i) \le -  \mu T/2 \right )  \\
& \le \frac{\sum_{i=0}^{T-1} \E \left ( \gD \hat Z(i) - \E  ( \gD \hat Z(i)|\fcal_i) \right)^2}{\mu^2 T^2/4 } = O(T^{-1})
\end{align*}

In the last line we use orthogonality of the increments, boundedness of $\gD \hat Z(t)$  and Chebyshev inequality.

Therefore we can conclude,
\begin{align*} \label{eq:dieprob}
\prob( X_1(T) + X_2(T)  \ge  \mu T/2 |\bfxx(0) =e_i ) & \ge \prob( Z(T) \ge  \mu T/2 |\bfxx(0) =e_i )  \\
& \ge  \prob(  \hat Z(T) \ge  \mu T/2,  \bfxx(t) \ne 0 \ \ \forall \ 0 \le t \le T | \bfxx(0) =e_i )\\
& \ge \prob(   \bfxx(t) \ne 0 \ \ \forall \ t \ge 0 | \bfxx(0) =e_i ) - \prob(  \hat Z(T) <  \mu T/2 | \bfxx(0) =e_i )\\
& \ge (1-q_i) - O(T^{-1}).
\end{align*}
The second inequality uses the fact that $Z(t)  = \hat Z(t)$ until  $\bfxx(t)$ reaches $(0,0)$.

\hfill $\square$

\section{Proof of the Theorem \ref{main:thm} part (b)}

\begin{lemma}
 In one round of exploration process on a subgraph involving $N  \ge (1-\gd/2)n$ many variables  the probability that (i) there is no termination due to the stopping rule and (ii)
 there  exists a contradictory cycle using the variables visited through the round is at least $p^2 \zeta^2$ for some $p, \zeta >0$.
 \end{lemma}

\noindent \textbf{Proof.} From Lemma \ref{dominate} and \ref{lem:sc}, we get
\begin{align*}
 \prob( a_T \ge C T, \tau =T)  & \ge \prob( X_1(T) + X_2(T) \ge C T , \tau =T)\\
  &\ge  \prob( X_1(T) + X_2(T) \ge C T)  -  \prob(  \tau <T)\\
  & \ge \eta - 2\exp(-\ga T/2)) \ge \zeta >0.
\end{align*}
  If  both $u$ and $\mom u \in A_T$ for some literal $u$, we have a directed path from the starting vertex to its complement using the literals in $E_T \cup A_T$. Otherwise, for each pair of literals $u, v \in A_T$, which are strongly distinct,  there are edges  $u \to \mom v$ and  $u \to \mom v$ in the digraph if the clause $(\mom u \vee \mom v)$ is present in the formula. If at least one of these $C T \choose 2$ many clauses is present, then  we again  have a directed path from  starting vertex to its complement using the literals in $E_T \cup A_T$.

\noindent {\bf Case I: $\mathbf {\ga_1 > 0} $.} Let $\ga_{\min} = \min( \ga_0, \ga_1, \ga_2) > 0$. Let $D$ be the event that $ \exists$  a directed path from  the starting vertex of the exploration process to its complement  in $E_T \cup A_T$. Then, by lemma \ref{lem:ind},

\begin{equation*}
 \prob( D  \big | \  a_T \ge CT  \text{ and }
 \tau= T ) \ge 1- (1 -\ga_{\min}/2n)^{CT \choose 2} \ge p >0.
 \end{equation*}

\noindent {\bf Case II: $\mathbf {\ga_1 = 0} $.}  Now $\ga_{\min} = 0$ and we can not make the above statement. But then, instead of looking at all $( \mom u \vee \mom v)$ clauses where $u,v$ are strongly distinct clauses belonging to $A_T$, we only consider such clauses where $u$ and $v$ have same parity. Since we have at least $2{CT/2 \choose 2}$ many clauses of that type, we have, similarly to case I,
\begin{equation*}
 \prob( D  \big | \  a_T \ge CT  \text{ and }
 \tau= T ) \ge  1 - (1 -\ga'/2n)^{ 2 {CT/2 \choose 2}} \ge p >0.
 \end{equation*}
where $\ga'= \min( \ga_0, \ga_2) > 0$.

Therefore,
\begin{equation*}
\prob( D, \tau=T )  \ge  \prob( D  \big | \  a_T \ge CT  \text{ and }
 \tau= T )   \prob( \  a_T \ge CT  \text{ and }
 \tau= T ) \ge p \zeta.\\
\end{equation*}

The lemma is now immediate from Corollary \ref{cor:ind}. \hfill $\square$

\begin{remark}
From the proof of the above lemma,  we have  seen that for large
$n$ with probability at least $r$ for some $r > 0$ we have a
directed path in the implication digraph from any literal $u$ to
its complement $\bar u$. Invoking FKG, we can say that we can find
a contradictory cycle with probability at least $r^2 >0$. But note
that we do not have a ready-made theorem like Friedgut's sharp
threshold result for generalized 2-SAT model.   Though we believe
that tweaking the lemma 4 of \cite{molloy} may help, we will not
pursue that. Instead we take a different route to bypass the
problem.
\end{remark}

\noindent  \textbf{Proof of the Theorem \ref{main:thm} part (b)} We now show how to bootstrap this positive probability event
 to an event with high probability.

 Initially we run a round of exploration process on the entire set of variables starting from $x_1$. If the process does not terminate after the first round, out of  at least $ n -4\alpha T -1$ many unvisited variables, we pick up an arbitrary one and run another round of exploration process in the deleted graph starting from it. We repeat this process $\gd \sqrt n/9\ga \le \gd n /2 (4 \ga T +1)$  many times, provided that we do not have to stop before, each time discarding previously visited variables to achieve independence among the different rounds. We thus ensure that in each run of exploration process, we have at least $( 1 - \gd/2)n$ many variables to start with.

We conclude that the probability that we get no contradictory cycle in all the rounds is at most
\begin{align*}
&\quad \  \ \prob( \text{no contradictory cycle and no round stops}) + \prob( \text{ one of the rounds stops})\\
&\le (1-p^2 \zeta^2)^{\gd \sqrt n/9\ga } +  (1- (1 - 2\exp(-\ga T/2))^{2\gd \sqrt n/9\ga } )\\
&\le \text{ Const} \times  \exp {(-B \sqrt n)}  \quad \text{for some  } B > 0.
\end{align*}
This concludes the proof. \hfill $\square$

\bibliographystyle{alpha}
\bibliography{sat}

\end{document}